\documentclass[graybox]{svmult}
\usepackage{fullpage}

\usepackage{type1cm}        
\usepackage{makeidx}         
\usepackage{graphicx}        
\usepackage{multicol}        
\usepackage[bottom]{footmisc}

\def\R{\mathbb{R}}
\def\N{\mathbb{N}}
\usepackage{url}
\usepackage{newtxtext}       %
\usepackage[varvw]{newtxmath}       
\usepackage{caption}
\makeindex             

\begin{document}

\title*{Nonlinear compressive reduced basis approximation for multi-parameter elliptic problem}
\titlerunning{Nonlinear compressive reduced basis method} 
\author{Hassan BALLOUT   \orcidID{0009-0005-5510-5672} \and\\ Yvon MADAY \orcidID{0000-0002-0443-6544} \and\\ Christophe PRUD'HOMME \orcidID{0000-0003-2287-2961}}
\institute{Hassan Ballout \at  Cemosis, IRMA UMR 7501, University of Strasbourg, CNRS, \email{hassan.ballout@etu.unistra.fr}
\and Yvon Maday \at Sorbonne Université, CNRS, Universit\'e Paris Cit\'e, Laboratoire Jacques-Louis Lions
(LJLL), \\F-75005 Paris, France \email{yvon.maday@sorbonne-universite.fr}
\and Christophe Prud'homme \at Cemosis, IRMA UMR 7501, Université de Strasbourg, CNRS, \email{prudhomme@unistra.fr}}
%
%
\maketitle

\abstract*{Blah}

\abstract  {: Reduced basis methods for approximating the solutions of parameter-dependant partial differential equations (PDEs) are based on learning the structure of the set of solutions - seen as a manifold ${\mathcal S}$ in some functional space - when the parameters vary. This involves investigating the manifold and, in particular, understanding whether it is close to a low-dimensional affine space. This leads to the notion of Kolmogorov $N$-width that consists of evaluating to which extent the best choice of a vectorial space of dimension $N$ approximates ${\mathcal S}$ well enough. If a good approximation of elements in ${\mathcal S}$ can be done with some well-chosen vectorial space of dimension $N$ -- provided $N$ is not too large -- then a ``reduced'' basis can be proposed that leads to a Galerkin type method for the approximation of any element in ${\mathcal S}$. In many cases, however, the Kolmogorov $N$-width is not so small, even if the parameter set lies in a space of small dimension yielding a manifold with small dimension. In terms of complexity reduction, this gap between the small dimension of the manifold and the large Kolmogorov $N$-width can be explained by the fact that the Kolmogorov $N$-width is linear while, in contrast, the dependency in the parameter is, most often, non-linear. There have been many contributions aiming at reconciling these two statements, either based on deterministic or AI approaches. We investigate here further a new paradigm that, in some sense, merges these two aspects:  the nonlinear compressive reduced basis
approximation. We focus on a simple multiparameter problem and illustrate rigorously that the complexity associated with the approximation of the solution to the parameter dependant PDE is  directly related to the number of parameters rather than the Kolmogorov $N$-width.}

\section{Introduction}
\label{sec:Intro}
From the very beginning,
the notion of ``learning'' is closely associated with --- if not at the root of ---  the reduced basis methods and  the complexity reduction techniques for the approximation of the solutions of parameter-dependant partial differential equations (PDEs) written as : 
\\
For some value of the parameter $\mu \in {\mathcal P}$ (the set of all parameters), find the solution $u_\mu\in X$ such that 
\begin{equation}\label{eq1}
    {\mathfrak R}(u_\mu;\mu)=0,
\end{equation} where $X$ is some appropriate Hilbert space, chosen for the above problem to be well-posed. 

In opposition to generic, multipurpose discretization methods for PDEs like finite volume, finite element, and spectral methods ... that can be used for any type of problem, the reduced basis methods are tuned to a certain class of problems of interest, the features of this problem are learned and digested, and the direct application to another problem makes no sense. These methods are based on a two phases approach : 
\begin{itemize}
    \item an ``{\sl off-line phase}'', that can be (and generally is) quite long and during which the learning of ${\mathcal S} =\{ u_\mu; \mu\in{\mathcal P}\}$ is performed
    \item an ``{\sl on-line phase}'', much (and generally very much) faster, during which the solution(s) $u_\mu$ for a new (set of) parameter(s) is (are) approximated.
\end{itemize}

One preliminary step in the --- {\sl off-line} --- learning process consists in investigating the set ${\mathcal S}$ of all solutions -- seen as a manifold in $X$ -- that consists in understanding whether the manifold is close to a low-dimensional affine space. This is related to the notion of Kolmogorov $N$-width that aims at evaluating the extent to which the best choice of an $N$-dimensional vector space $X_N$ approximates ${\mathcal S}$ sufficiently well. The discovery of 
such an accurate $N$-dimensional vector space $X_N$ is not straightforward and getting the best one may even be impossible. Various algorithms (POD, SVD, Greedy \cite{bookRb, bookCrb}),  have been proposed to determine alternative quasi-optimal $X_N$, and performed in the off-line phase based on the data of many instances in ${\mathcal S}$.
Then, provided that  such a quasi-optimal $X_N$ has been identified,  a ``reduced'' basis method can be deployed --- in the {\sl on-line phase} ---  under the form, e.g. of
a Galerkin formulation in $X_N$ for the approximation of the solution $u_\mu$ for some given value $\mu$ of the parameter. This low complexity approximation works for a large variety of parameter-dependant problems and provides fast, even, superfast approximation methods, with a complexity that scales as ${\mathcal O}(N^3)$. 

 In many instances however, even if the set of parameters lies in a space of small dimension ${\mathcal P} \subset \R^P$ ($P\in \N$),  the Kolmogorov $N$-width of ${\mathcal S}$ is not so small (even if ${\mathcal S}$ can then be considered as a manifold with dimension $\le P$). This conceptual mismatch, between the manifold dimension and the complexity Kolmogorov width can be understood by the fact that the Kolmogorov width is a linear concept whereas the dependency of the solution $u_\mu$ in the parameter $\mu$ is, most often highly non-linear. However, it is true that the intrinsic complexity of the manifold is, heuristically at least, related to the dimension of ${\mathcal P}$.  
 
 Several authors have proposed nonlinear approximation/compression approaches to deal with this problem. There are essentially 2 wide classes of approaches: Lagrangian and Eulerian. 

Lagrangian approaches (also called registration methods see \cite{Taddei2020}) involve a map $\varphi_\mu$ so that the set of all $v_\mu = u_\mu \circ \varphi_\mu$ defines a manifold better suited to linear compression methods than the original ${\mathcal S}$ with thus a faster decay of the  Kolmogorov width. Examples of Lagrangian approaches have been proposed in \cite{iollo2014, mojgani2017, cagniart2019, ohlberger2013, Taddei2015, ferrero2022}. Eulerian approaches have a more nonlinear root: some approaches are based on a partition of the parameter domain in an adaptive way \cite{eftang2011}, others rely on basis updates \cite{maday2013, carlberg2015, etter2019, peherstorfer2020}, or on a Grasmannian learning \cite{amsallem2008, zimmermann2018, polack2021} or finally based on artificial learning methods, some involving convolutional autoencoders, \cite{kashima2016, Youngkyu2020, Lee2020, Fresca2021} or Operator Inference \cite{Kramer2024}.

In these last references, the original notion of ``learning'' of the reduced basis methods is hence pushed further to the extent of investigating the nonlinearity of the manifold. This learning process can be used to the point of not even referring to the problem \eqref{eq1} and at the end propose a method that, to any input parameter $\mu\in {\mathcal P}$, proposes -- as an output -- an approximation to $u_\mu$. Even if this is shown to work in those publications, this approach has some obvious drawbacks : 
 \begin{itemize}
     \item the size of the database for the learning stage, needs to be quite large in order to apprehend the complexity of $\mathcal{S}$
     \item there is no evaluation of the residual of equation  \eqref{eq1} that would allow to state 
     \begin{itemize}
         \item if the problem is well or not approximated 
         \item what to do in order to improve a non-sufficient accuracy
     \end{itemize}
     
 \end{itemize}
 this is why, following \cite{barnett2022quadratic, Geelen2023}, the authors in \cite{barnett2023, cohen2023nonlinear} have proposed to ``help'' the statistical approaches by using intermediate models and still ground the approximation of the approximated resolution of problem \eqref{eq1}.
 
In this paper, we investigate further this idea applied to a simple multiparameter problem and show, in particular, that, indeed, the complexity associated with the approximation of the solution to the parameter-dependent PDE is related directly to the number of parameters rather than the Kolmogorov $N$-width.

\bigskip

The paper is organized as follows: Section 2 presents the problem setting, detailing the multiparameter elliptic problem and the physical model used for our analysis. Section 3 provides a comprehensive discussion on the analysis of the solution manifold, introducing the classical reduced basis methods and extending the discussion to our proposed nonlinear compressive approach. Section 4 explores the application of machine learning techniques for nonlinear approximation within the reduced basis framework, including detailed discussions on the regression models used for the recovery of higher RB modes from lower ones. Section 5 illustrates the efficacy of the proposed method through a series of numerical experiments, demonstrating the advantages of our approach in terms of accuracy and computational efficiency. Finally, Section 6 concludes the paper with a summary of our findings and suggestions for future research directions. 

The numerical experiments presented in the paper are conducted using the Feel++ library~\cite{feelpp} to produce the finite element method (FEM) computations as well as the reduced basis using proper orthogonal decomposition (POD). For the machine learning aspects of our methodology, we utilize the scikit-learn framework~\cite{sklearn}.

\section{Problem setting: the multiparameter problem}
\label{sec:2}
\begin{figure}[htbp]
\centering
\includegraphics[width=0.5\textwidth]{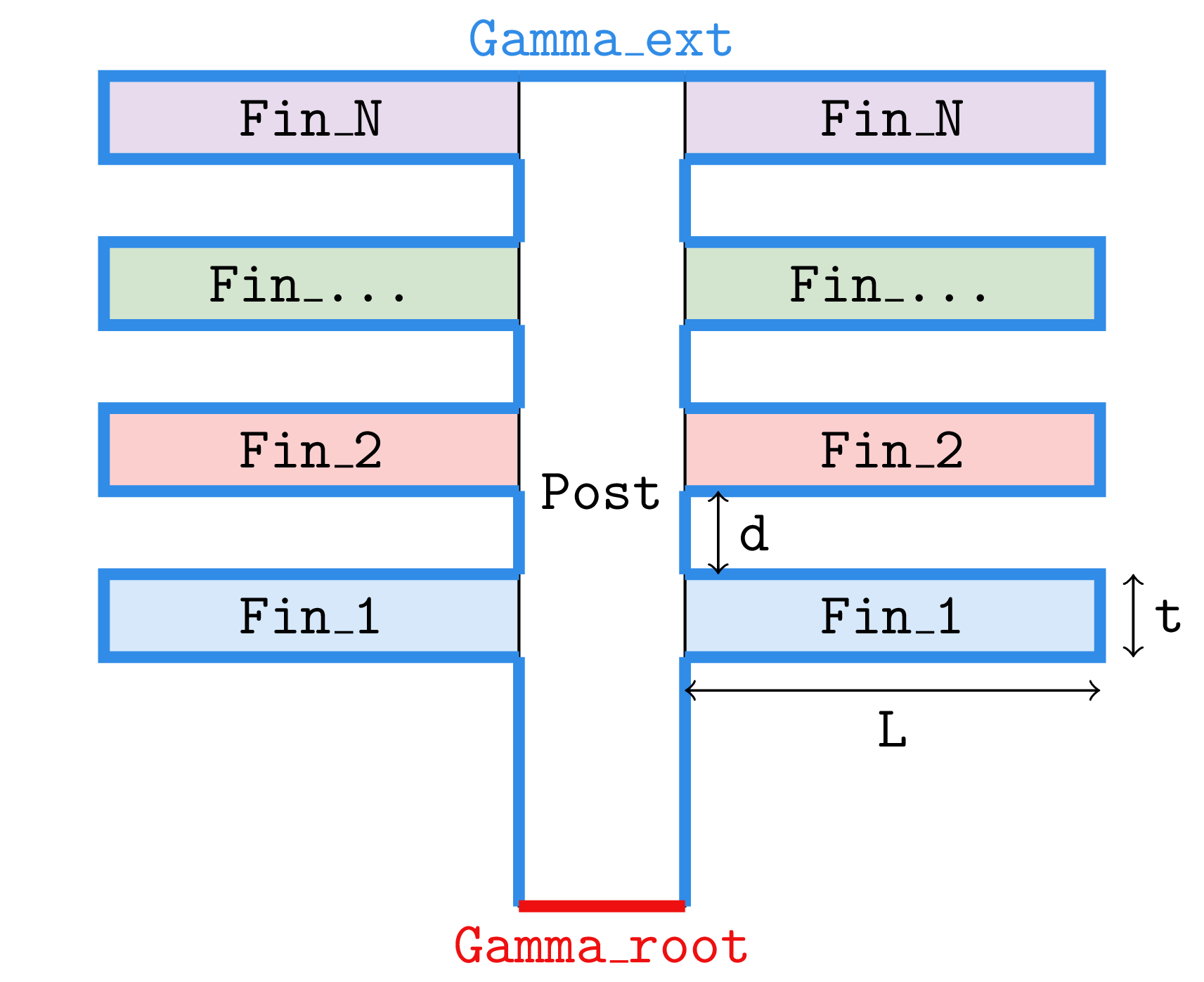}
\caption{Schematic of a thermal fin showing the structure and important elements, see also \cite{rbpp}.}
\label{fig:thermal-fin}
\end{figure}

We consider the problem of designing a thermal fin to remove heat from a surface effectively, see {\emph e.g.} \cite{rbpp}. The two-dimensional fin, shown in Figure~\ref{fig:thermal-fin}, consists of a vertical central \textbf{post} and say $N_{fins}$ horizontal \textbf{subfins}; the fin conducts heat from a prescribed uniform flux source at the root, $\mathbf{\Gamma_{\text{root}}}$, through the large-surface-area subfins to surrounding flowing air. The fin is characterized by a ($P=N_{fins}+1$)-component parameter vector $\mu = (\mu_1, \mu_2, \ldots, \mu_{P})$, where $\mu_i = k^i$, $i = 1, \ldots, N_{fins}$, and $\mu_{P} = Bi$; $\mu$ may take on any value in a specified design set $\mathcal{P} \subset \mathbb{R}^{P}$.

Here $k^i$ is the thermal conductivity of the $i$-th subfin (normalized relative to the post conductivity $k^0 = 1$); and $Bi$ is the Biot number, a nondimensional heat transfer coefficient reflecting convective transport to the air at the fin surfaces (larger $Bi$ means better heat transfer). For example, suppose we choose a thermal fin with $N_{fins}=4$, $k^1 = 0.4$, $k^2 = 0.6$, $k^3 = 0.8$, $k^4 = 1.2$, and $Bi = 0.1$; for this particular configuration $\mu = (0.4, 0.6, 0.8, 1.2, 0.1)$, which corresponds to a single point in the set of all possible configurations $\mathcal{P}$. 

We are interested in the design of this thermal fin, and we thus need to look at certain outputs or cost-functionals of the temperature as a function of $\mu$. For our output $T_{\text{root}}$, we choose the average steady-state temperature of the fin root normalized by the prescribed heat flux into the fin root. The particular output chosen relates directly to the cooling efficiency of the fin --- lower values of $T_{\text{root}}$ imply better thermal performance. The steady-state temperature distribution within the fin, $u_\mu$, is governed by the elliptic partial differential equation
\begin{equation}
\label{eq:testcase-strong-Omega}
- \nabla . [{\bf k}\nabla] u_\mu = 0 \  \text{ in } \Omega, \, 
\end{equation}
which implies in particular that
\begin{equation}
\label{eq:testcase-strong}
-k^i\Delta u = 0 \text{ in } \Omega^i, \, i = 0, \ldots, N_{fins},
\end{equation}
where $\Delta$ is the Laplacian operator, and $u_i$ refers to the restriction of $u$ to $\Omega^i$. Here $\Omega^i$ is the region of the fin with conductivity $k^i$, $i = 0, \ldots, N_{fins}$; $\Omega^0$ thus the central post, and $\Omega^i$, $i = 1, \ldots, N_{fins}$, corresponds to the subfins. The entire fin domain is denoted $\Omega = \bigcup_{i=0}^{N_{fins}} \Omega^i$; the boundary is denoted $\mathbf{\Gamma}$. The continuity of temperature and heat flux at the conductivity-discontinuity interfaces is also implicitly set on $\Gamma_{int}^i \equiv \partial\Omega^0 \cap \partial\Omega^i$, $i = 1, \ldots, N_{fins}$, where $\partial\Omega^i$ denotes the boundary of $\Omega^i$, we have on $\Gamma^i$, $i = 1, \ldots, N_{fins}$:
\begin{equation}
\label{eq:interface-conditions}
\begin{aligned}
[u] &= 0, \\
\left[k^i \frac{\partial u}{\partial n^i}\right] &= 0,
\end{aligned}
\end{equation}
where $n^i$ is the outward normal on $\partial\Omega^i$. Finally, we introduce a Neumann flux boundary condition on the fin root:
\begin{equation}
\label{eq:flux-bc}
-k^0 \frac{\partial u}{\partial n^0} = -1 \text{ on } \Gamma_{\text{root}}. 
\end{equation}
which models the heat source; and a Robin boundary condition:
\begin{equation}
\label{eq:robin-bc}
-k^i \frac{\partial u}{\partial n^i} = Bi\,u^i  \text{ on } \Gamma^i_{\text{ext}}, \, i = 1, \ldots, N_{fins},
\end{equation}
which models the convective heat losses. Here $\Gamma^i_{\text{ext}}$ is that part of the boundary of $\Omega^i$ exposed to the flowing fluid; note that $\bigcup_{i=0}^{N_{fins}}\Gamma^i_{\text{ext}} = \Gamma\setminus\Gamma_{\text{root}}$. The average temperature at the root, $T_{\text{root}}(\mu)$, can then be expressed as $\ell^O(u_\mu)$, where
\begin{equation}
\label{eq:testcase-output}
\ell^O(v) = \int_{\Gamma_{\text{root}}} v.
\end{equation}
Note that $\ell(v) = \ell^O(v)$ for this problem (compliant problem).

In order to proceed with a Galerkin type approximation of the solutions to problem \eqref{eq:testcase-strong}, it is now written under a  parametric variational problem :
Given $\mu \in {\mathcal P}$, find $u_\mu \in X$ such that:
\begin{equation}
  \label{eq:varpb}
 a(u_\mu,v;\mu) = f(v;\mu), \quad \forall v \in X
\end{equation}
here $a(.,.;\mu)$ is thus a symmetric elliptic bilinear form.

\section{Analysis of the solution manifold}
\label{sec:3}
\subsection{Kolmogorov width - Classical RBM approximation}
\label{subsec:3.1}

The  Kolmogorov $m$-width of ${\mathcal S}$ is defined as
\begin{equation}\label{kw}
    d_m({\mathcal S}) = 
    \inf_{\dim(X_m)\leq m}  \max_{v\in {\mathcal S}}
    \min_{w\in X_m} \| v -w\|_X.
\end{equation}
As was recalled in the introduction, this describes how well the manifold can be approximated
by an ideally selected (and usually out of reach) $m$-dimensional space. A slow decay of the Kolmogorov $m$-width means that a large number of basis functions is needed to achieve a good approximation, which compromises the efficiency of our reduced basis method.

The determination of the optimal spaces that realize the minimum is most of the times an open problem and has been circumvented, from a practical point of view, within the frame of classical reduced basis methods, particularly Proper Orthogonal Decomposition (POD) that aim to approximate the solution manifold by a well-chosen (quasi-optimal) N-dimensional vector space.

For the sake of completeness, let us recall the definition of the POD (or SVD) method. We start by introducing a  training set ${\mathcal P}_h = \{\mu_1, \mu_2, ..., \mu_M\} \subset {\mathcal P}$ of $M$ parameter values, with $M$ chosen large enough to ``represent well enough'' the set of parameters ${\mathcal P}$. 

For each parameter value $\mu_i \in {\mathcal P}_h$, we approximate the corresponding solution $u_{\mu_i}\in {\mathcal S}$ with, -- say -- an accurate enough finite element method (hence the presence of the index $h$), that we denote as $u_{h,\mu_i}$.
Then, we construct the correlation matrix $C \in \mathbb{R}^{M \times M}$ as follows:

\begin{equation}
  \label{eq:correlationmatrix}
  C_{ij} = \langle u_{h,\mu_i}, u_{h,\mu_j} \rangle_{X} \quad \forall i,j \in \{1,2,...,M\}
\end{equation}
where $\langle .,.\rangle_{X}$ denotes the inner product in $X$.
We consider the eigenvalue-eigenvector pairs $(\lambda_i, V_i)$ of the correlation matrix $C$ such that:
\begin{equation}
  \label{eq:eigenvalueproblem}
  C \, V_i = \lambda_i \, V_i, \quad \forall i \in \{1,2,...,M\}
\end{equation}
and:
\begin{equation}
  \label{eq:eigenvalueorder}
  \lambda_1 \geq \lambda_2 \geq ... \geq \lambda_M
\end{equation}
The POD basis functions $\zeta_i$ corresponding to the eigenvectors  $V_i$ are given by:

\begin{equation}
  \label{eq:PODbasis}
  \zeta_i = \frac{1}{\lambda_i} \sum_{j=1}^{M} (V_{i})_j \, u_{h,\mu_j} \quad \forall i \in \{1,2,...,M\}
\end{equation}
where $(V_{i})_j$ denotes the $j$-th component of $V_i$.

The key property of this basis is that if one truncates the POD basis to the first $N$ functions $\{\zeta_1,\zeta_2, ..., \zeta_N\}$, they span the optimal $N$-dimensional subspace $X_N$. This optimality is in the sense of minimizing the following error measure:

\begin{equation}
  \label{eq:PODerror}
   \sum\limits_{\mu \in {\mathcal P}_h} \inf\limits_{{v_N}\in X_N} \|u_{h,\mu} - v_N\|^2_X
\end{equation}
over all $N$-dimensional subspaces $X_N$ of $X_M= \text{span}\{u_{h,\mu} \, | \, \mu \in {\mathcal P}_h\}$.
In practice, we choose only $N \ll M$ modes corresponding to the largest N eigenvalues $\lambda_i$. This reveals a drawback of the POD method: it requires solving $M\gg N$ truth problems in order to explore the solution manifold. An alternative iterative method is the Greedy method.

For our multiparameter problem, with 8 subfins $(N_{fins}=8)$ and a piecewise linear finite element functions space of dimension of $N_h \approx 10000$. We checked the quality of the finite element solution previously described for $N_{\mathrm{fins}}=8 (P=9)$ for the Kolmogorov width estimation. 
The average $L_2$ error norm of the finite element solutions described previously with respect to a finer mesh and quadratic finite element instead of linear is about $3. 10^{-3}$ in average with a maximum value of $5. 10^{-3}$ over the sample of parameters used.

\begin{figure}[htbp]
    \centering
\includegraphics[width=0.7\textwidth] {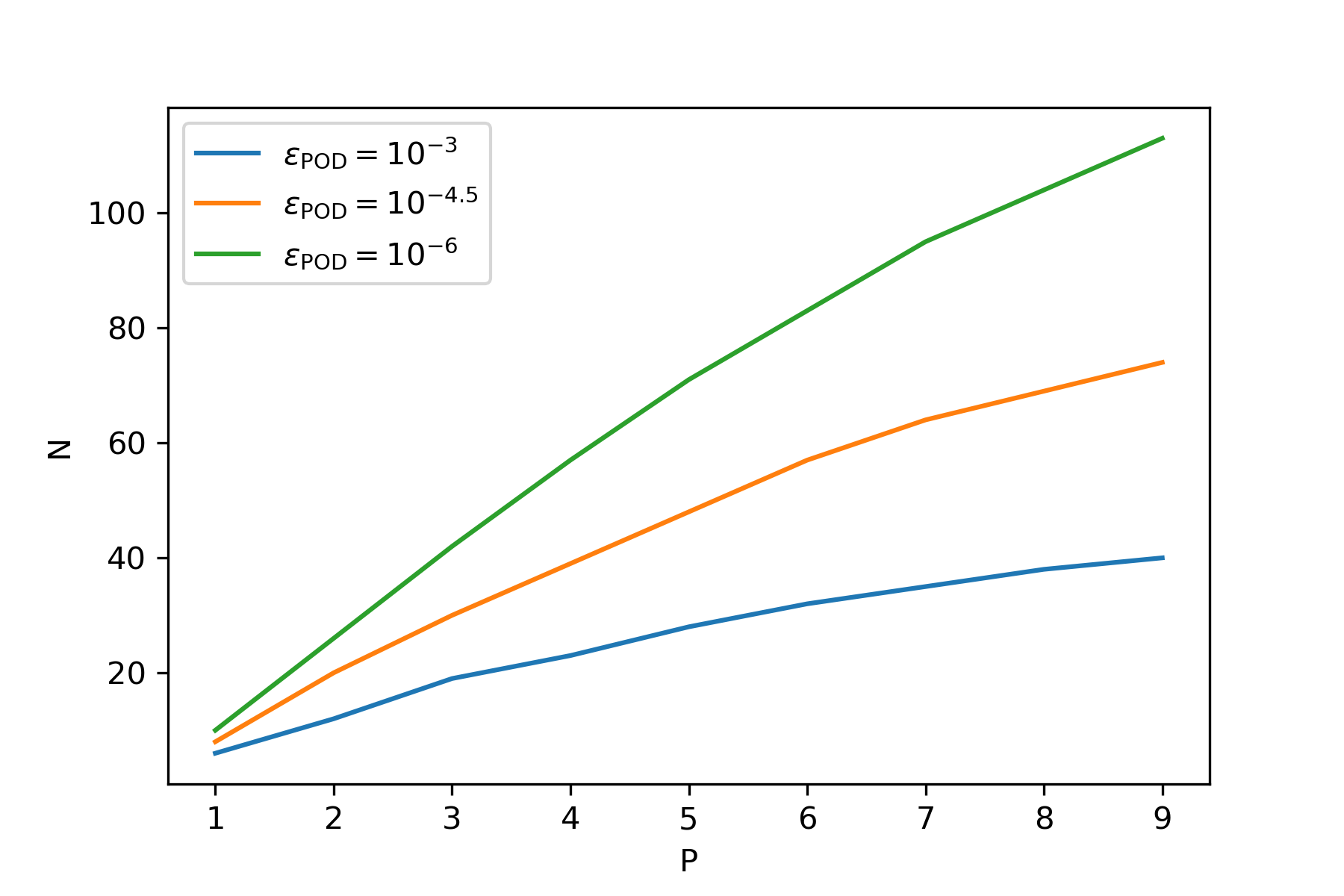}
    \caption{Number of reduced basis functions as a function of the parameter dimension for different accuracies}
    \label{fig:NvsP}
\end{figure}

We show in Figure \ref{fig:NvsP} the number of POD modes, $N$, needed to reach a given accuracy, \(\epsilon_{POD}\), as a function of the parameter space dimension, $P$. We can see that the Kolmogorov $N$-width highly depends on $P$. More precisely, as $P$ increases, the decay of the Kolmogorov $N$-width becomes slower (even if it can be described as fast). Note that the notion of $\epsilon_{POD}$ is based on the Relative Information Content, and $N(\epsilon)$ is chosen as the smallest value that satisfies:
$$ I(N) = \frac{\sum_{i=1}^N \lambda_i}{\sum_{i=1}^M \lambda_i} \geq 1 - \epsilon^2 .$$

Following this purely POD-based analysis to estimate the Kolmogorov width, we can achieve a reduced basis approximation using a POD-Galerkin method. 

We thus defined the reduced space $X_N = \text{span}\{\zeta_1,\dots,\zeta_N\} \subset X$. Given a parameter value $\mu$ in ${\mathcal P}$,  we seek a reduced basis approximation $u^N_\mu$ in $X_N$ such that:
\begin{equation}
  \label{eq:reducedvarpb}
  a(u^N_\mu,v_N;\mu) = f(v_N;\mu), \quad \forall v_N \in X_N
\end{equation}
This reduced basis approximation $u^N_\mu$ is expressed as follows:
\begin{equation}
  \label{eq:reducedapprox}
  u^N_\mu = \sum_{i=1}^{N} \hat u_{i,\mu} \zeta_i
\end{equation}
And by testing the previous formulation with $v_N=\zeta_i$ for $i=1,\dots,N$, we obtain the coordinates vector of $u^N_\mu$ as the solution of the following linear system:
\begin{equation}
  \label{eq:reducedmatrix}
  \mathbf{A}_{N}^{\mu}\mathbf{u}_{N}^{\mu} = \mathbf{f}_{N}^{\mu}
\end{equation}
where $\mathbf{A}_{N}^{\mu}=\{a(\zeta_i,\zeta_j;\mu)\}_{i,j=1}^{N}$ and $\mathbf{f}_{N}^{\mu} = \{f(\zeta_i;\mu)\}_{i=1}^{N}$. Finally, we can evaluate the reduced basis output as:
\begin{equation}
  \label{eq:reducedoutput}
  s_N(\mu) = (\mathbf{u}_{N}^{\mu})^T \mathbf{f}_{N}^{\mu}
\end{equation}

\begin{figure}[htbp]
    \centering
\includegraphics[width=1.0\textwidth] {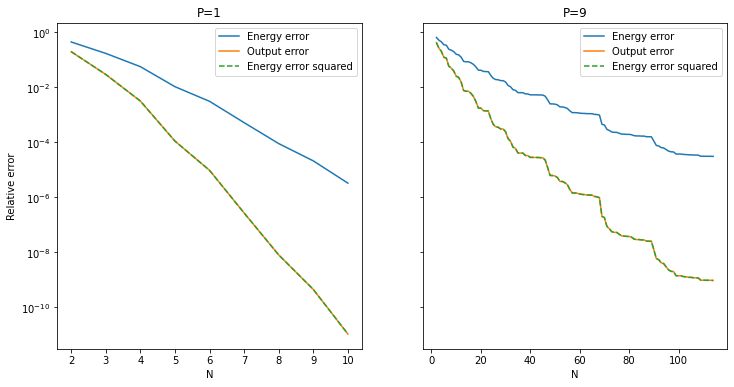}
    \caption{Mean relative errors vs number of POD modes for $P = 1$(left) and $P = 9$(right)}
    \label{fig:P1P9errorsVsN}
\end{figure}

In Figure \ref{fig:P1P9errorsVsN}, we illustrate the decay of both energy and output errors with respect to $N$ for $P = 1$ and $P = 9$. By reducing the dimension from $N_h$ to $N$, we obtained a speed-up factor of 1000 for $P = 1 \,(N = 10)$ with a mean relative energy error (between the reduced method approximation and the finite element solution) of $10^{-6}$, and a speed-up factor of 100 for $P = 9$ $(N = 114$) with a relative energy error of $10^{-5}$.  Additionally, we observe that the output error aligns with the square of the energy error, which is a classic result for compliant problems.

From now on, we use a coarser grid with linear piecewise finite element and $N_{\mathrm{fins}}=5.$ 
The average $L_2$ error norm with respect to the fine grid setting with quadratic piecewise finite elements is then in average $7e-3$ and the maximum value $1e-2$ over the sample of parameters used.

The reduced basis space is constructed using POD such that $\epsilon_{POD} = 10^{-6}$. Note that, of course, the overall error with respect to the exact solution is, in our case, dominated by the finite element error. In contrast, the subsequent errors in our nonlinear compressive reduced basis are computed with respect to the finite element approximation.

\subsection{Linear and nonlinear notions of complexity reduction}
\label{subsec:3.2}
In order to investigate the nonlinear complexity reduction, let us remind that the frame of complexity reduction is classically described by a pair of continuous mappings, the encoder
$$
E: X \to \R^m,
$$
and the decoder
$$
D: \R^m \to X.
$$
These two mappings are tuned so that the quantity
\begin{equation}\label{eqED}
    \max_{v\in{\mathcal S}} \|v-D(E(v))\|_X \hskip0.1truecm, \hskip0.5truecm  \hbox{is minimized under some constraints,}
\end{equation}
 and $m$ here aims at reflecting the complexity  of ${\mathcal S}$. 
 
 For instance, the already mentioned Kolmogorov width $d_m({\mathcal S})_X$ assumes that both $E$ and $D$ are linear mappings.

Even if these two notions should not be mixed up, in order to reconcile the dimension $P$ of the manifold ${\mathcal S}$ and the complexity reduction it is natural to think at alternative definitions of the complexity (see e.g. \cite{devoreActa98, devoreActa21} and also \cite{cohen2009} for the link with Gelfand width) that leads to the idea of enlarging the class where the encoder and decoder stay (see \cite{cohen2023nonlinear} for an extended presentation in this frame). This is somehow the idea explored in references \cite{Lee2020, Fresca2021}, where the notion of auto-encoder is used to minimize the information needed to reconstruct the identity in ${\mathcal S}$, as a combination of two neural networks ${\mathcal E}$ and ${\mathcal D}$ (i.e. $Id = {\mathcal D} \circ {\mathcal E}$), and where the size of the layer at the output of ${\mathcal E}$ and the input of ${\mathcal D}$ is minimized.. looking thus to identify $m$. In these publications, neural networks are therefore used to discover ${\mathcal E}$ and ${\mathcal D}$, and further investments are made in the field of machine learning, with two warnings: the amount of data required to train the associated networks, and the lack of understanding underlying machine learning.

In the opposite direction, we can refer to a, still timely but much older problem of ``parameter reconstruction'' or ``parameter identification'', and ask what minimal knowledge on the elements $v$ of ${\mathcal S}$ is required to identify the value of parameters $\mu$ in  ${\mathcal P}$ such that $v=u_\mu$. Indeed, assuming that we have such a reconstruction  ${\mathcal R} : v\in {\mathcal S}\mapsto \mu\in {\mathcal P}$, then the identity mapping in ${\mathcal S}$ : $Id_{\mathcal S}$ would involve the solution operator $\mu\in {\mathcal P}\mapsto u_\mu\in {\mathcal S}$ so that $Id_{\mathcal S} = v\mapsto u_{{\mathcal R}(v)}$.

Our approach stands in the middle of these two extreme encoder-decoder paradigms and relies on the two notions.
The first one is the {\it sensing numbers} $s_m({\mathcal S})$
where $E$ in \eqref{eqED}  is still assumed to be linear but $D$ is assumed
to be nonlinear. In other words, they can 
be defined as

\begin{equation}
    \label{eqsensing}
    s_m({\mathcal S}):=
\inf_{D,\ell_1,\dots,\ell_m}
\max_{v\in{\mathcal S}} \|v-D(\ell_1(v),\dots,\ell_m(v))\|_X,
\end{equation}
where the infimum is taken over all choices of linear functionals $\ell_1,\dots,\ell_m\in X'$
and decoding map $D$. Then, referring to the parameter reconstruction/identification, we can raise the question: what are the best linear mappings $\ell_1,\dots,\ell_m$ that represent well ${\mathcal S}$, or, if we remember that we are in a Hilbert setting, from Riesz identification between $X$ and $X'$, the question becomes: what are the best elements in $X$: $v_1,\dots,v_m$ that 
are the more aligned with ${\mathcal S}$ (in order to have the largest scalar product with any element of ${\mathcal S}$), so that, for any $w\in {\mathcal S}$ and any $j, 1\le j\le m$: $\ell_j(w)= \langle v_j,w\rangle_{X}$ are the best linear functional appearing above? This is actually what is considered in the POD framework within the notion of Relative Information Content: the first modes on the POD approach are the ones that contain in the vectorial space spanned by ${\mathcal S}$ the largest amount or information. This naturally leads to considering that the first POD functions $v_j=\zeta_j$ are the best elements that carry the information on solutions in $ {\mathcal S}$. This is what led us in \cite {barnett2022quadratic, barnett2023, cohen2023nonlinear} to propose the ``anszatz'' that a candidate for the optimal encoder in \eqref{eqsensing} consists in the moments associated with those POD modes.

In order to evaluate this claim, let us understand to which extent -- and how many of them -- these first POD moments allow to recover the parameters and can define a reconstruction ${\mathcal R} : v\in {\mathcal S}\mapsto \mu\in {\mathcal P}$, that is written under the form ${\mathcal R} = {\mathcal D}\circ E$
where the encoder $E$ is linear and defined as being the set of moments
\begin{equation}
    v\mapsto \{ \langle \zeta_i,v\rangle_{X}\}_{i=1,..,m}
\end{equation}

Actually, in dimension $P=1$ (only the Biot number is varying) the following figure \ref{fig:dataP1Inverse}
graphically illustrates the simple dependence of the Biot number as a function of the first moment of the solution $u_{Bi}$, which is the sought for Decoder ${\mathcal D}$. The (Lipschitz) continuity of $u$ with respect to $Bi$ and, thus, the continuity of the decoder $D$, is a classical result (see e.g. \cite{bookRb} Sect.5.3 for the proof).

\begin{figure}[!htbp]
    \centering
\includegraphics[width=0.7\textwidth] {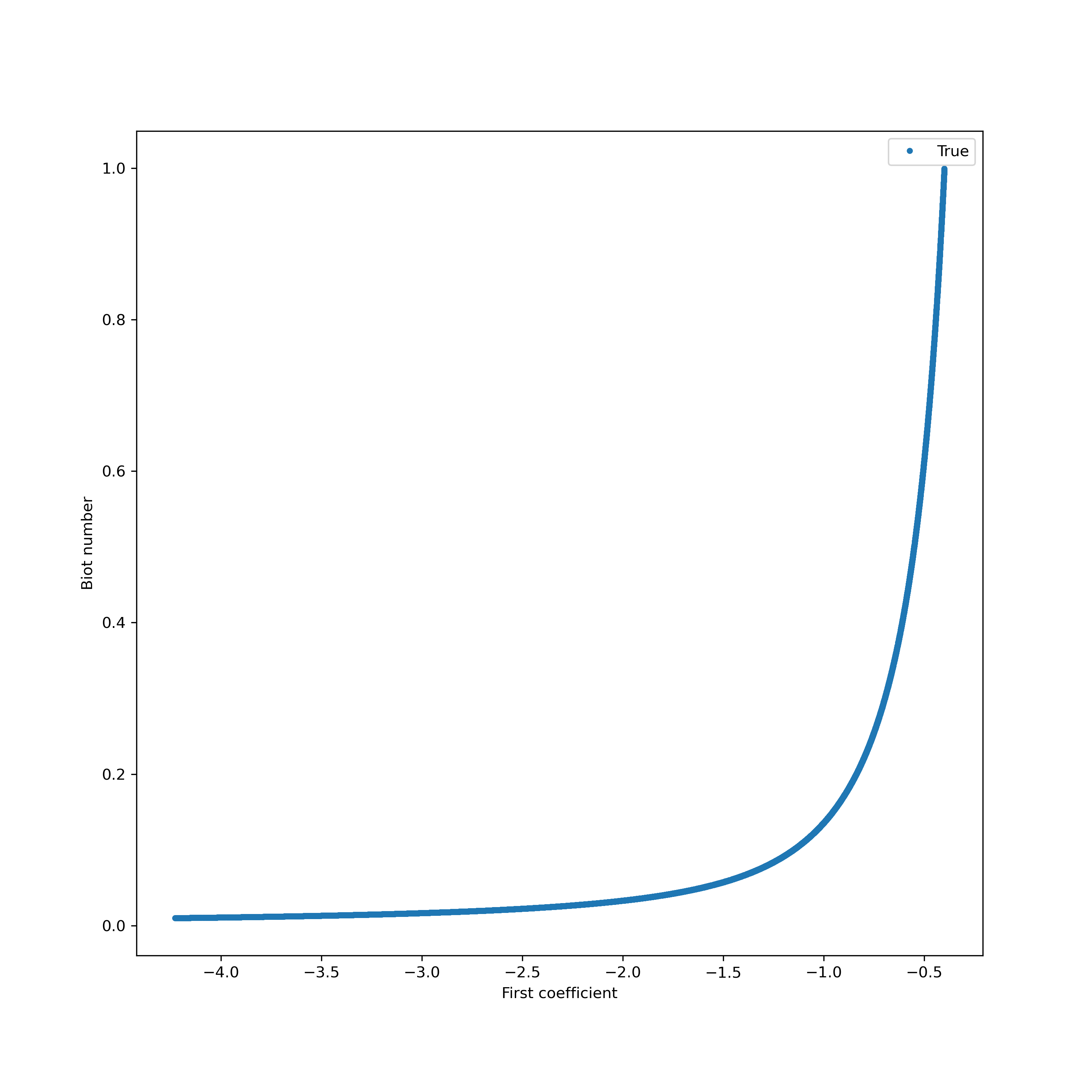}
    \caption{ The Biot number as a function of the first RB coefficient}
    \label{fig:dataP1Inverse}
\end{figure}

\subsection{Sensing number and Parameter identification}
\label{subsec:3.3}

This being said, in higher dimensions, we have to investigate the nonlinear decoder. As is presented in \cite {barnett2022quadratic, barnett2023, cohen2023nonlinear}, where polynomial and ensemble learning regression methods (like a random forest) can be used to define the decoder.

Back with the 1-dimensional problem ($P=1$),  the dependence of $\mu = Bi$ with respect to the first coefficient displayed 
in Figure \ref{fig:dataP1Inverse}
 can very well be approximated as a spline regression problem. By increasing the number of knots ($n_{\text{knots}}$) and the degree of the model, spline regression can achieve machine precision. Using a training dataset size of $10^7$, decision trees and random forests also deliver excellent accuracy, with mean relative errors reaching up to $10^{-7}$. However, polynomial models perform very weakly on this problem and do not achieve the same level of precision.

In higher-dimensional problems ($P>1$), the regression problem becomes harder, and some algorithms suffer from the curse of dimensionality, as is well reported in the literature
see~\cite{curse}. In our case, the reconstruction of the parameter using decision trees (or random forests) gives "acceptable" results up to $P=3$ (a mean relative error of $10^{-3}$ on the Biot number reconstruction and $10^{-2}$ for $K^0$ and $K^1$). For the same reason, for $P>3$, more advanced models need to be employed to achieve a reasonable reconstruction of the parameters.

We can conclude from these initial results that indeed $s_m({\mathcal S}) = 0$ for $m \geq P$, at least for $\ P=1,2,3$, which sets the connection between the manifold dimension and the complexity of the manifold written in terms of sensing number rather than Kolmogorov $m$ width.

To finish, we have also numerically investigated the stability of the inverse problem for $P=1,2,3$ for the abovementioned regressors. The inverse problem is unstable whenever we use the polynomial regressors and their input data is outside the training domain; otherwise, in all other cases, it is stable. The inverse problems associated with spline, decision tree, and random forest regressors are, on the contrary, stable, including for input data outside the training domain. 

\section{Recovering the high RB modes from the low ones}
\label{sec:4}
After this analysis on the inverse/identification problem, let us come back to the approximation of the functions in ${\mathcal S}$. One of the fundamental steps in the nonlinear compressive reduced basis method is the regression task, specifically training a model that accurately predicts the high RB coefficients from the first $n$ independent ones. Both accuracy and performance (prediction cost for approximating the high modes from the low ones) are crucial in this step. Accuracy is vital because any error in this step introduces a corresponding error in the final RB solution. Performance is 
equally important, as it significantly impacts online performance, especially if the prediction is used iteratively. Another fundamental consideration is the stability of these approximations. An unstable model is more likely to cause divergence in the associated solvers.

\subsection{$P=1$}
\label{subsec:4.1}
Starting with the 1-dimensional case ($P=1$), we consider the problem of predicting the last $N-n$ coefficients from the first $n$ coefficient(s). Following the discussion in sections 3.2 and 3.3, it is natural to expect that $n=1$ is sufficient for the reconstruction of high modes. The value of $N$ is determined such that $\epsilon_{POD} = 10^{-6}$, which gave us $N(\epsilon) = 9$.

\begin{figure}[!hb]
    \centering
\includegraphics[width=0.9\textwidth] {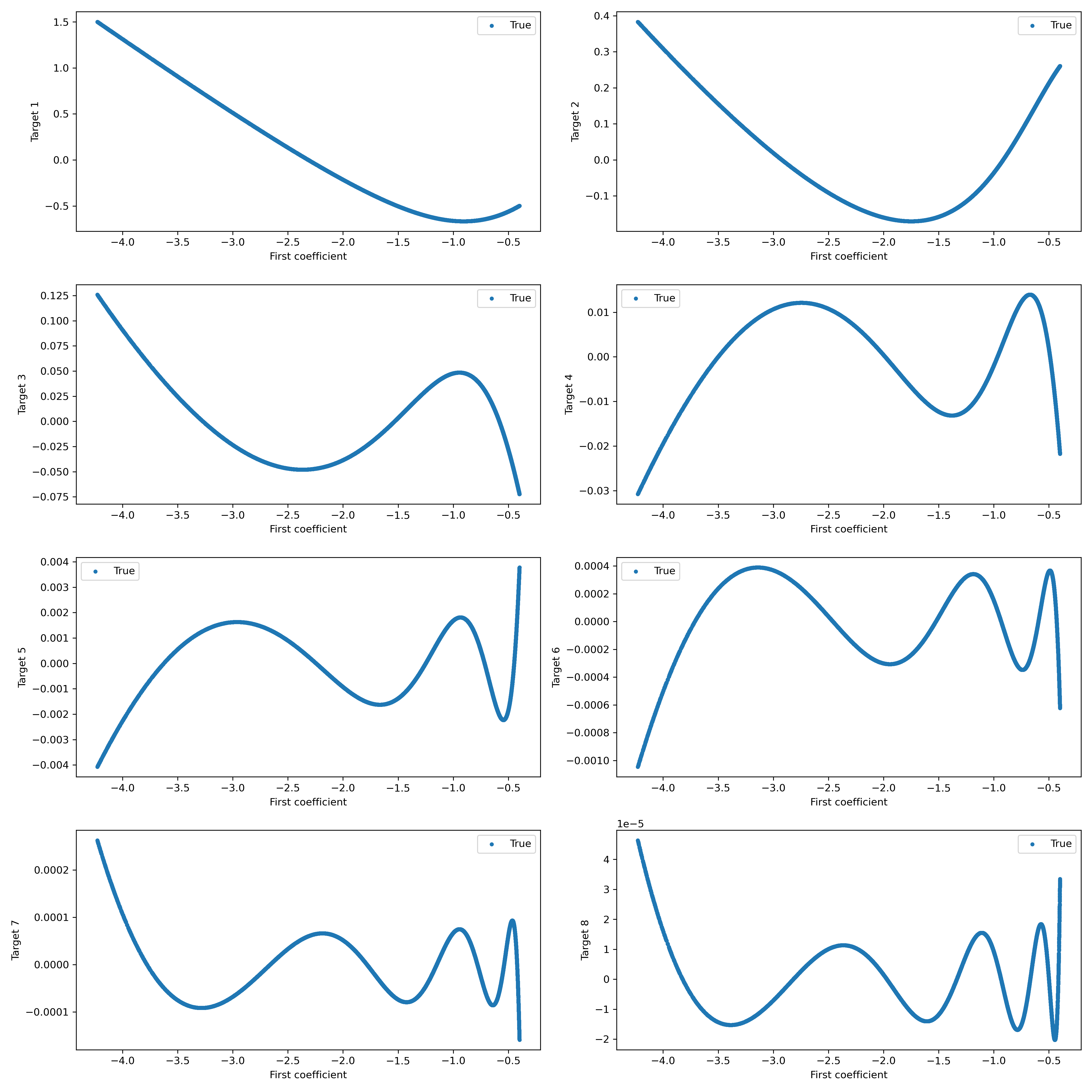}
    \caption{ Higher coefficients as a function of the first coefficient for $P=1$ and $n=1$}.
    \label{fig:dataP1}
\end{figure}

\begin{figure}[!htbp]
    \centering
\includegraphics[width=0.7\textwidth] {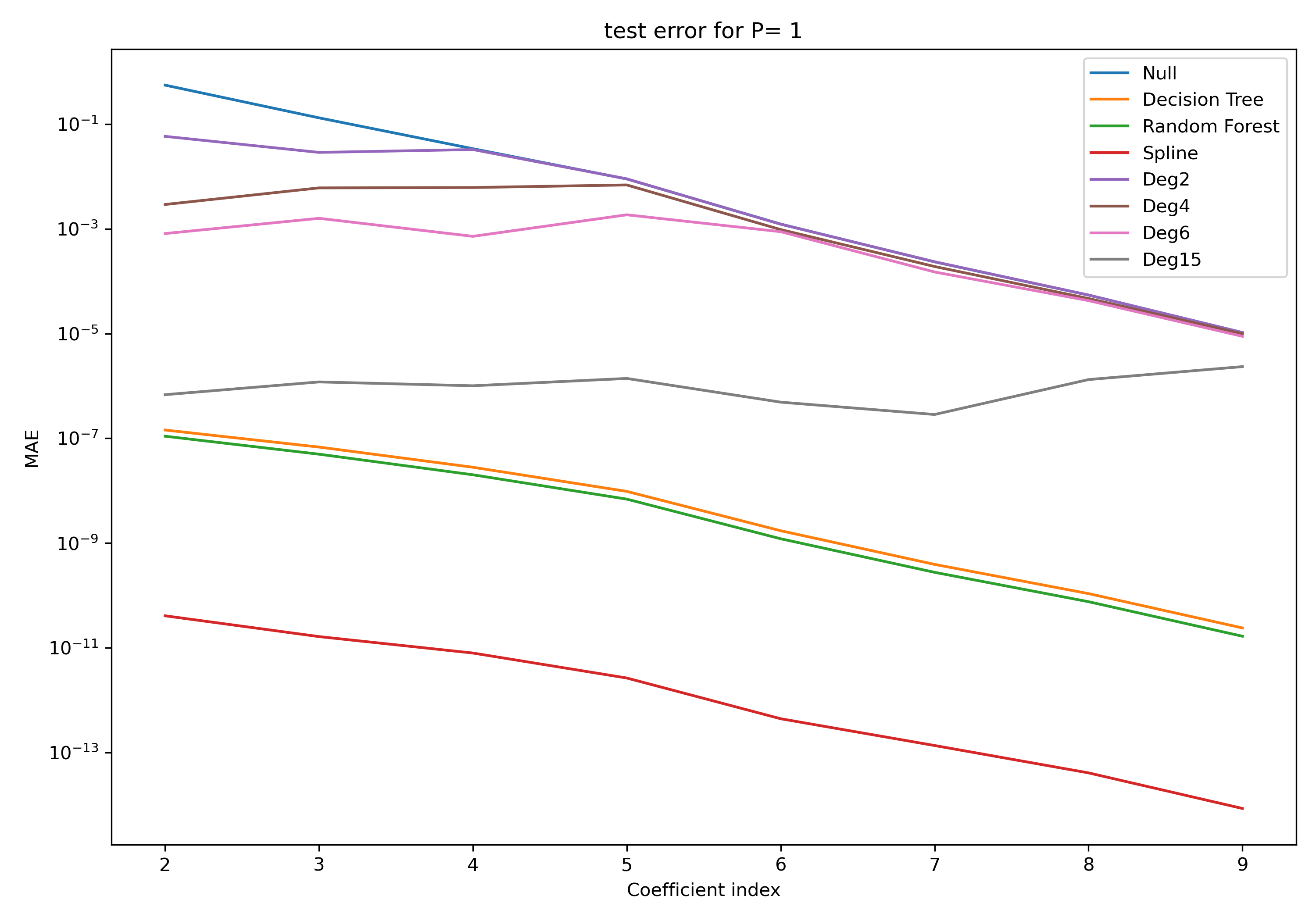}
    \caption{Mean absolute error on a test set for different ML models for $P=1$ and $n=1$.}
    \label{fig:MLP1}
\end{figure}

For this regression task, we tested several classical ML algorithms, namely decision tree, random forest, spline, and polynomial regression. The training dataset size is $10^7$ for decision tree and random forest models and $10^5$ for polynomial and spline models. In Figure \ref{fig:MLP1}, we can see that splines provide the best accuracy compared to other models. The accuracy of polynomial regression improves with increasing degree due to the oscillatory behavior of the higher RB coefficients (see Figure \ref{fig:dataP1}). The spline model, tuned with hyperparameters (knots=100, degree=10), achieves the best accuracy due to its local adaptability, giving it a significant advantage over polynomials.

%

\begin{figure}[!hb]
    \centering
\includegraphics[width=1.0\textwidth] {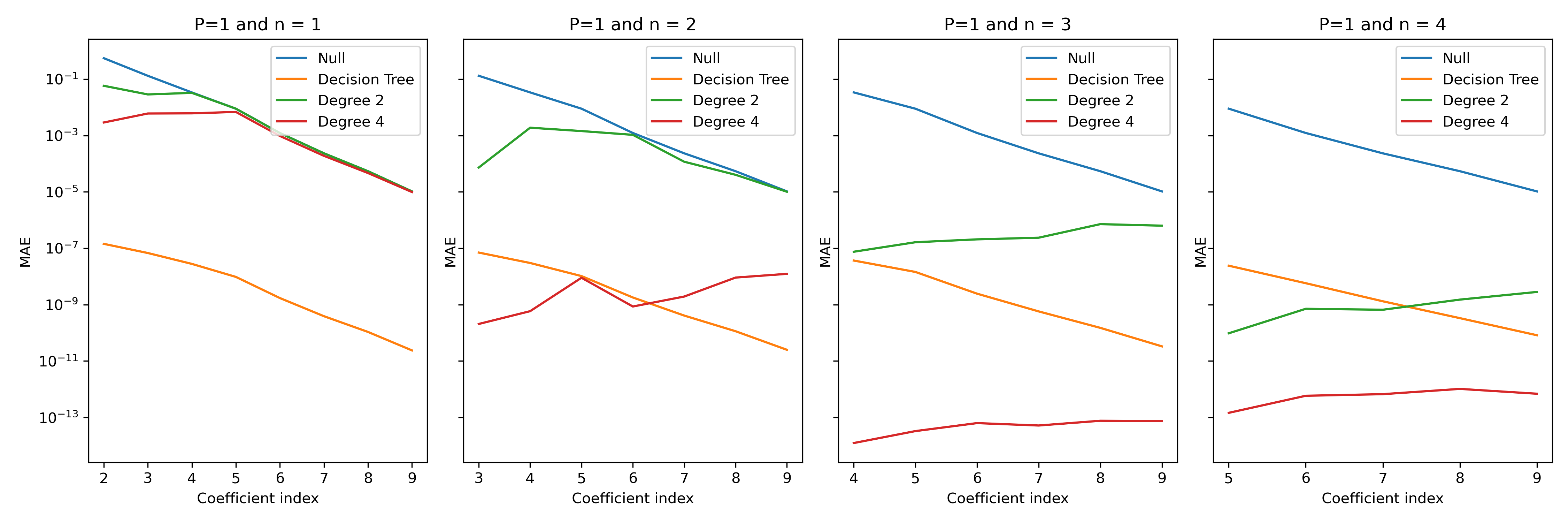}
    \caption{Mean absolute error for different number of features $n$ and $P=1$}
    \label{fig:P1n}
\end{figure}

Next, we consider increasing values of $n$ for $P=1$. The results are depicted in Figure~\ref{fig:P1n}. We observe that decision trees (and similarly random forests) do not improve, while polynomial regressors gain in accuracy. This result is expected and can be explained by the fact that the first coefficient is the only independent one. Increasing $n$ does not provide additional information but implicitly increases the degree of polynomial models, which explains the improvement in accuracy.

\begin{figure}[!hb]
    \centering
\includegraphics[width=1.0\textwidth] {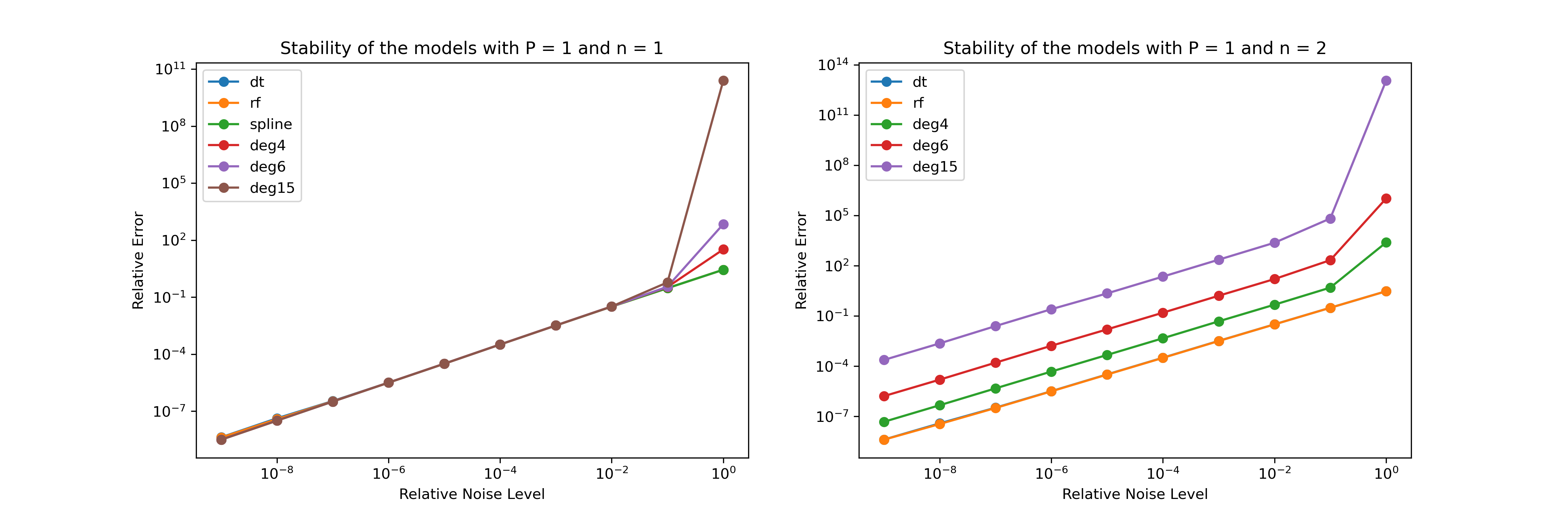}
    \caption{Stability of models for $P=1$, $n=1$ (left panel) and $n=2$ (right panel).}
    \label{fig:stabP1}
\end{figure}

While this approach appears to work well with polynomials, it has significant drawbacks. Firstly, increasing $n$ significantly raises the prediction cost, given by $\binom{n+d}{d}$ where $d$ is the degree of the polynomial. Secondly, it introduces substantial instability in the model (see Figure~\ref{fig:stabP1}), making it impractical for the online phase where a nonlinear system is solved iteratively. 

\subsection{$P>1$}
\label{subsec:4.2}
\begin{figure}[!hb]
    \centering
\includegraphics[width=0.7\textwidth] {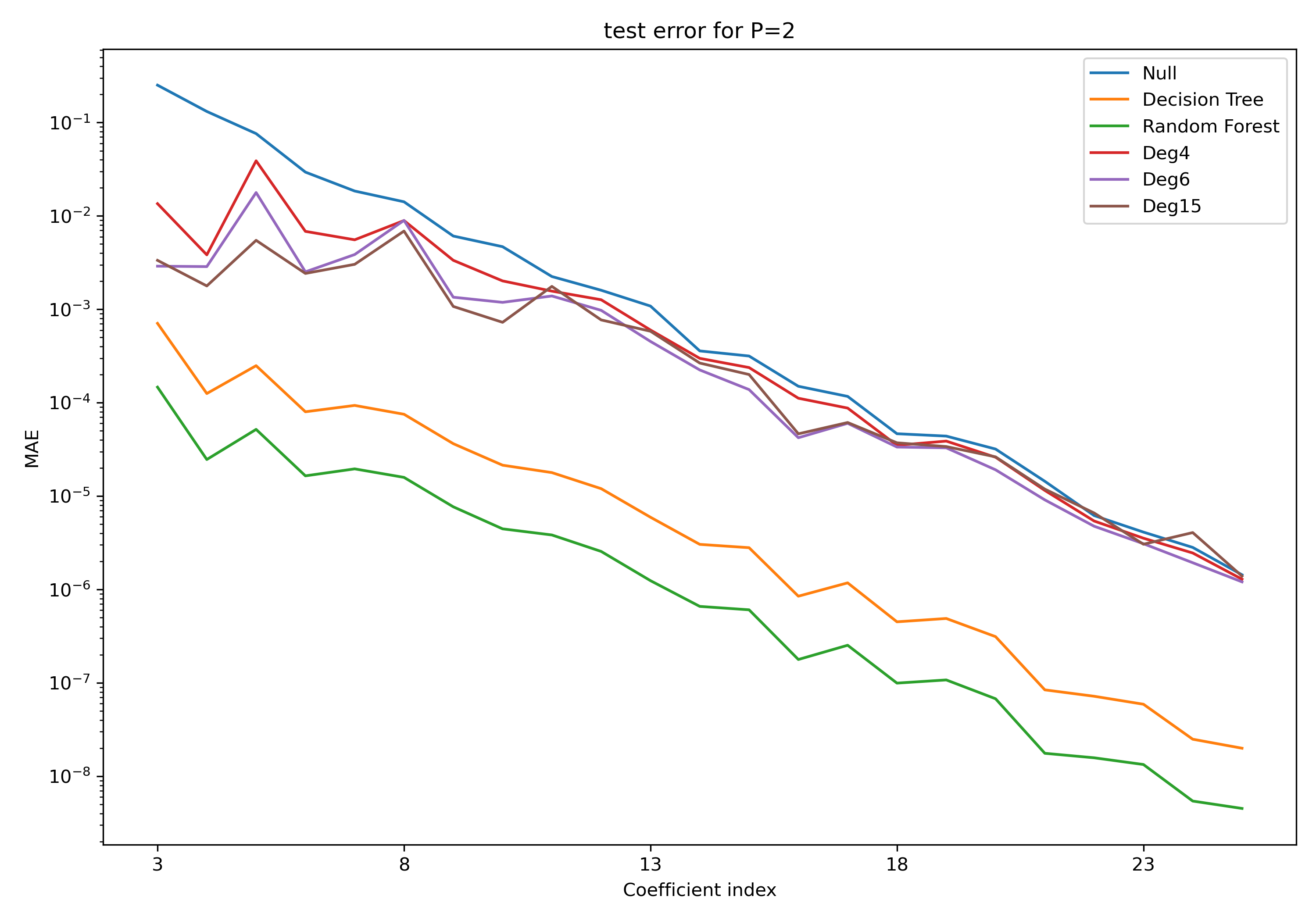}
   \caption{Mean absolute error on a test set for different ML models for $P=2$}
    \label{fig:MLP2}
\end{figure}

We now consider the case $P=2$. The number of POD modes $N(\epsilon_{POD})=25$ where $\epsilon_{POD}=10^{-6}$. Figure~\ref{fig:MLP2} shows that for $n=P=2$, decision trees and random forests perform better than polynomials; however, we do not achieve the same accuracy as for $P=1$. Similar to $P=1$, increasing $n$ improves the accuracy of polynomial models but leads to instability.

For higher parameter dimensions ($P>2$), as seen in section \ref{subsec:3.1}, the multivariate regression task becomes more complex. Achieving acceptable accuracy using these classic machine-learning algorithms is challenging and improved approaches need be involved.

\section{Nonlinear compressive RBM in action}
\label{sec:5}
Once we have identified proper (linear) encoder and (nonlinear) decoder  to approximate elements in ${\mathcal S}$, we can use this paradigm to approximate the solution of problem \eqref{eq1} for new values of $\mu$. This is done by writing, e.g. : Find $u^{n,N}_\mu$
 in $X_N$ that is written as
\begin{equation}\label{eqCRBsol}
   u^{n,N}_\mu = \sum_{i=1}^n \hat u_{i,\mu} \zeta_i +  \sum_{j=n+1}^N  \hat u_{j}(\{\hat u_{i,\mu}\}_{i=1,..,n})
\zeta_j
\end{equation} such that
\begin{equation}\label{eqCRBeq}
  \Pi_n\Bigl[ {\mathfrak R}(u^{n,N}_\mu;\mu)\Bigr]=0,
\end{equation} 
where $\Pi_n$ is some rank $n$ operator. 

In this equation, the degress of freedom are the $n$ values $\{\hat u_{i,\mu}\}_{i=1,..,n}$, and the unknown $\hat u_{j}(\{\hat u_{i,\mu}\}_{i=1,..,n})$ for $j=n+1,\dots, N$ are functions of the degrees of freedom. In our case, we choose a Galerkin approach and choose $ \Pi_n$ to be the projection on $X_n$.

This problem is nonlinear and iterative methods are advocated to solve it.

The problem is quite simple, and we have implemented a Picard iteration procedure (for a larger problem, a Newton or quasi-Newton method should be implemented). This gives rise to the following (where $k$ is the iteration index and where we leave all dependencies in $\mu, n$ and $N$).

\begin{equation}
  \sum_{i=1}^n \hat u_{i}^{k+1} \zeta_i  =
  \sum_{i=1}^n \hat u_{i}^{k} \zeta_i
  - \gamma^k \Pi_n\Bigl[ {\mathfrak R}(u^k;\mu)\Bigr]
\end{equation}
with
\begin{equation}
     u^k = \sum_{i=1}^n \hat u_{i}^k \zeta_i +  \sum_{j=n+1}^N  \hat u_{j}(\{\hat u_{i}^k\}_{i=1,..,n})
\zeta_j
\end{equation}

To present the results, in Section \ref{subsec:4.1}, we showed that for $P=1$, the best regression accuracy is achieved using the spline model. Applying the Nonlinear Compressive Method with $n=1$ yields the results shown in Table~\ref{tab:errorsP1}. We observe that the nonlinear reconstruction of high modes significantly reduces the error compared to the classical RB method with only $n=1$ mode. It is also evident that regression accuracy plays a crucial role in the final results, and with splines, we can achieve nearly the same accuracy as RB with $N=9$ modes.

\begin{table}[h]
    \centering
    \begin{tabular}{|l|c|c|}
        \hline
        \textbf{Model} & \textbf{Relative Energy Error} & \textbf{Relative Output Error} \\
        \hline
        RB (with 9 modes) & $9.39 \times 10^{-7}$ & $1.50 \times 10^{-12}$ \\
        \hline
        RB (with 1 mode) & $4.68 \times 10^{-1}$ & $2.89 \times 10^{-1}$ \\
        \hline
        Degree 4 regression & $1.17 \times 10^{-2}$ & $1.31 \times 10^{-3}$ \\
        \hline
        Degree 6 regression & $2.91 \times 10^{-3}$ & $6.48 \times 10^{-4}$ \\
        \hline
        Degree 15 regression & $3.41 \times 10^{-6}$ & $4.70 \times 10^{-7}$ \\
        \hline
        Spline & $9.39 \times 10^{-7}$ & $2.24 \times 10^{-11}$ \\
        \hline
        Decision tree & $9.91 \times 10^{-7}$ & $4.72 \times 10^{-8}$ \\
        \hline
        Random forest & $9.78 \times 10^{-7}$ & $3.70 \times 10^{-8}$ \\
        \hline
    \end{tabular}
    \captionsetup{justification=centering}
    \caption{Relative energy and output Errors for P=1}
    \label{tab:errorsP1}
\end{table}

For $P=2$, a mean energy error of $10^{-4}$ can be achieved using a decision tree model, compared to $10^{-1}$ with 2 modes and $10^{-6}$ with 25 modes using the classical RBM.  In this case and for $P>2$, it is essential to employ more robust nonlinear solvers alongside advanced regression models to achieve better accuracy.

We want finally and consider the potential economy from the pure reduced basis method to the nonlinear compressive reduced basis approximation.

The plain reduced basis method involves :
\begin{itemize}
    \item as already indicated in the introduction, for a linear problem: the inversion of a full $N\times N$ matrice, i.e. ${\mathcal O}(N^3)$ operations
    \item for nonlinear problems, thanks to hypereductions methods like EIM \cite{EIM2004}, only an addition (and negligible) contribution like ${\mathcal O}(N^2)$ 
\end{itemize}
On the other hand, this nonlinear compressive RBM involves
\begin{itemize}
    \item for a linear problem : the inversion of a $n\times n$ matrix hence ${\mathcal O}(n^3)$ operations and for each iteration $k=1,..,K$ the evaluation of the high modes with respect to the $n$ first ones, hence ${\mathcal O}(n^2 + (N-n)^2)$ complexity
    \item and additional $(K {\mathcal O}( N^2))$, again if the original problem is nonlinear.
\end{itemize}
which of course rapidly came in favor of this latter approach.

\section{Conclusion and perspective}
\label{sec:Conclusion}

As a preamble to this conclusion, we would like to point out that the proposed approach can be considered from two angles. The first, and the one we prefer, is that of approximating elements $u_\mu$ of ${\mathcal S}$: given some data that are linear evaluations of $u_\mu \in {\mathcal S}$ (the first $n$ coefficients of $u_\mu$ in a reduced basis), how can we construct a good approximation of $u_\mu$ by authorizing non-linear constructions (of the next $N-n$ coefficients of $u_\mu$ in the reduced basis) from the data? We're well within the “sensing numbers” framework and this approximation paradigm is used to determine the first $n$ coefficients of $u_\mu$. 

The other way of looking at the approach is to think in terms of $\mu$ and not $u_\mu$. We then look directly at the problem \eqref{eqCRBeq}, which is a nonlinear problem in the degrees of freedom $\{\hat u_{i,\mu}\}_{i=1,..,n}$ and the linearity of the encoder is hidden. We want to stress that the situation is similar if we consider a classical RBM approximation for a nonlinear problem. the resulting system is nonlinear in the component of the solution, while the foundation of the approach is properly linear and based on Kolmogorov's $N$-width.

Going back to the purpose of this paper, we have demonstrated on the low dimension $P=1,2,3$ at least, that we can indeed compute through a non-linear construction $N-n$ coefficients from the first $n$ coefficients associated with the more important modes. The critical offline ingredient is the regressor fitted to reconstruct the $N-n$ coefficients online with the Picard iterations. A few remarks are necessary. Several regressors have been used, --- decision trees and random forests, polynomials are various degrees and splines --- and compared. As $P$ increases, the regressor accuracy decreases, or the regressor is not available, such as the spline which gave the best results for $P=1$. Clearly, more advanced regressors must be used to support higher dimensions. The online procedure will also need to be studied in more detail concerning the robustness and computational cost. Finally, comparing the present work with other generation strategies for the reduced basis space would be interesting, {\emph e.g.} greedy.

\begin{acknowledgement}
This work has been partially funded by the European Research Council (ERC) under the European
Union's Horizon 2020 research and innovation program (grant No 810367), project EMC2 (YM), and the France 2030 NumPEx Exa-MA (ANR-22-EXNU-0002) project managed by the French National Research Agency (ANR).

The authors want to thank Gérard Biau, Anthony Nouy and Agustin Somacal for interesting discussions during the preparation of this manuscript.
\end{acknowledgement}

\end{document}